\input amstex
\documentstyle {amsppt}
\pagewidth{12.5 cm}\pageheight{19 cm}\magnification\magstep1
\topmatter
\title Unipotent variety in the group compactification \endtitle
\author Xuhua He\endauthor
\affil Department of Mathematics, M.I.T., Cambridge, MA 02139
\endaffil \keywords Unipotent varieties; Steinberg fibers; Group
compactifications \endkeywords \subjclassyear{2000} \subjclass
20G15 \endsubjclass
\thanks {\it E-mail address}: hugo\@math.mit.edu \endthanks

\abstract The unipotent variety of a reductive algebraic group $G$
plays an important role in the representation theory. In this
paper, we will consider the closure $\bar{\Cal U}$ of the
unipotent variety in the De Concini-Procesi compactification
$\bar{G}$ of a connected simple algebraic group $G$. We will prove
that $\bar{\Cal U}-\Cal U$ is a union of some $G$-stable pieces
introduced by Lusztig in \cite{L4}. This was first conjectured by
Lusztig. We will also give an explicit description of $\bar{\Cal
U}$. It turns out that similar results hold for the closure of any
Steinberg fiber in $\bar{G}$.
\endabstract

\endtopmatter
\document

\define\po{\text{\rm pos}}

\define\supp{\text{\rm supp}}
\define\End{\text{\rm End}}

\define\Lie{\text{\rm Lie}}
\define\Spec{\text{\rm Spec}}

\define\Ad{\text{\rm Ad}}
\redefine\i{^{-1}}
\redefine\ge{\geqslant}
\redefine\le{\leqslant}

\define\cb{\Cal B}

\define\cn{\Cal N}

\define\cp{\Cal P}

\define\cu{\Cal U}

\define\a{\alpha}
\redefine\b{\beta}

\define\g{\gamma}

\redefine\d{\delta}
\redefine\D{\Delta}
\define\e{\epsilon}

\redefine\o{\omega}
\define\p{\pi}

\define\s{\sigma}

\redefine\l{\lambda}

\redefine\v{\vartheta}

\define\tz{\tilde Z}
\define\tl{\tilde L}
\define\tc{\tilde C}

\head Introduction \endhead

A connected simple algebraic group $G$ has a ``wonderful''
compactification $\bar{G}$, introduced by De Concini and Procesi.
The variety $\bar{G}$ is a smooth, projective variety with $G
\times G$ action on it. The $G \times G$-orbits of $\bar{G}$ are
indexed by the subsets of the simple roots.

The group $G$ acts diagonally on $\bar{G}$. Lusztig introduced a
partition of $\bar{G}$ into finitely many $G$-stable pieces. The
$G$-orbits on each piece are in one-to-one correspondence to the
conjugacy classes of a certain reductive group. Based on the
partition, he developed the theory of ``Parabolic Character
Sheaves'' on $\bar{G}$.

In this paper, we study the closure $\bar{\cu}$ of the unipotent
variety $\cu$ of $G$ in $\bar{G}$, partially based on the previous
work of \cite{Spr2}. The main result is that the boundary of the
closure is a union of some $G$-stable pieces. (see Theorem 4.3.)

The unipotent variety plays an important role in the
representation theory. One would expect that $\bar{\cu}$, the
subvariety of $\bar{G}$, which is analogous to the subvariety
$\cu$ of $G$, also plays an important role in the theory of
``Parabolic Character Sheaves''. Our result is a step toward this
direction.

The arrangement of this paper is as follows. In section 1, we
briefly recall some results on the $B \times B$-orbits of
$\bar{G}$ (where $B$ is a Borel subgroup of $G$) and results on
$\bar{\cu}$, which were proved by Springer in \cite{Spr1} and
\cite{Spr2}. In section 2, we first recall the definition of the
$G$-stable pieces and then in 2.6, we show that any $G$-stable
piece is the minimal $G$-stable subset of $\bar{G}$ that contains
a particular $B \times B$-orbit. In the remaining part of section
2, we establish some basic facts about the Coxeter elements, which
will be used in section 4 to prove our main theorem. In section 3,
we show case-by-case that certain $G$-stable pieces are contained
in $\bar{\cu}$. Hence a lower bound of $\bar{\cu}$ is established.

A naive thought about $\bar{\cu}$ is that the boundary of the
``unipotent elements'' are ``nilpotent cone''. In fact, it is
true. A precise statement is given and proved in 4.3. Thus we
obtain an upper bound of $\bar{\cu}$. We also show in 4.3 that the
lower bound is actually equal to the upper bound. Therefore, our
main theorem is proved. In section 4, we also consider the closure
of arbitrary Steinberg fiber of $G$ in $\bar{G}$. (An example of
Steinberg fiber is $\cu$.) The results are similar. In the end of
section 4, we calculate the number of points of $\bar{\cu}$ over a
finite field. The formula bears some resemblance to the formula
for $\bar{G}$.

\head 1. Preliminaries \endhead

\subhead 1.1 \endsubhead Let $G$ be a connected, simple algebraic
group over an algebraically closed field $k$. Let $B$ be a Borel
subgroup of $G$, $B^-$ be the opposite Borel subgroup and $T=B
\cap B^-$. Let $(\a_i)_{i \in I}$ be the set of simple roots. For
$i \in I$, we denote by $\a^\vee_i$, $\o_i$, $\o^\vee_i$ and $s_i$
the simple coroot, the fundamental weight, the fundamental
coweight and the simple reflection corresponding to $\a_i$. We
denote by $< , >$ the standard pairing between the weight lattice
and the root lattice. For any element $w$ in the Weyl group
$W=N(T)/T$, we will choose a representative $\dot w$ in $N(T)$ in
the same way as in \cite{L1, 1.1}.

For any subset $J$ of $I$, let $W_J$ be the subgroup of $W$
generated by $\{s_j \mid j \in J\}$ and $W^J$ (resp. $^J W$) be
the set of minimal length coset representatives of $W/W_J$ (resp.
$W_J \backslash W$). Let $w^J_0$ be the unique element of maximal
length in $W_J$. (We will simply write $w^I_0$ as $w_0$.) For $J,
K \subset I$, we write $^J W^K$ for $^J W \cap W^K$.

\subhead 1.2 \endsubhead For $J \subset I$, let $P_J \supset B$ be
the standard parabolic subgroup defined by $J$ and $P^-_J \supset
B^-$ be the opposite of $P_J$. Set $L_J=P_J \cap P^-_J$. Then
$L_J$ is a Levi subgroup of $P_J$ and $P^-_J$. Let $Z_J$ be the
center of $L_J$ and $G_J=L_J/Z_J$ be its adjoint group. We denote
by $\p_{P_J}$ (resp. $\p_{P^-_J}$) the projection of $P_J$ (resp.
$P^-_J$) onto $G_J$.

Let $\bar{G}$ be the wonderful compactification of $G$ (\cite{DP}
deals with the case $k=\bold C$. The generalization to arbitrary
$k$ was given in \cite{Str}). It is an irreducible, projective
smooth $G \times G$-variety. The $G \times G$-orbits $Z_J$ of
$\bar{G}$ are indexed by the subsets $J$ of $I$. Moreover, $Z_J=(G
\times G) \times_{P^-_J \times P_J} G_J$, where $P^-_J \times P_J$
acts on the right on $G \times G$ and on the left on $G_J$ by $(q,
p) \cdot z=\p_{P^-_J}(q) z \p_{P_J}(p) \i$. Let $h_J$ be the image
of $(1, 1, 1)$ in $Z_J$.

We will identify $Z_I$ with $G$ and the $G \times G$-action on it
is given by $(g, h) \cdot x=g x h \i$.

For any subvariety $X$ of $\bar{G}$, we denote by $\bar{X}$ the
closure of $X$ in $\bar{G}$.

For any finite set $A$, we will write $|A|$ for the cardinality of
$A$.

\subhead 1.3 \endsubhead For any closed subgroup $H$ of $G$, we
denote by $H_{diag}$ the image of the diagonal embedding of $H$ in
$G \times G$ and by $\Lie(H)$ the corresponding Lie subalgebra of
$H$. For $g \in G$, we write $^g H$ for $g H g \i$.

For any parabolic subgroup $P$, we denote by $U_P$ its unipotent
radical. We will simply write $U$ for $U_B$ and $U^-$ for
$U_{B^-}$. For $J \subset I$, set $U_J=U \cap L_J$ and $U^-_J=U^-
\cap L_J$.

For parabolic subgroups $P$ and $Q$, define $$P^Q=(P \cap Q)
U_P.$$

It is easy to see that for $J, K \subset I$ and $u \in {}^J W^K$,
$P_J^{(^{\dot u} P_K)}=P_{J \cap \Ad(u) K}$.

Let $\cu$ be the unipotent variety of $G$. Then $\cu$ is stable
under the action of $G_{diag}$ and $U$ is stable under the action
of $U \times U$ and $T_{diag}$. Moreover, $\cu=G_{diag} \cdot U$.
Similarly, $\bar{\cu}=G_{diag} \cdot \bar{U}$ (see \cite{Spr2,
1.4}).

\subhead 1.4 \endsubhead Now consider the $B \times B$-orbits on
$\bar{G}$. We use the same notation as in \cite{Spr1}. For any $J
\subset I$, $u, v \in W$, set $[J, u, v]=(B \times B) (\dot{u},
\dot{v}) \cdot h_J$. It is easy to see that $[J, u, v]=[J, x, v z
\i]$, where $u=x z$ with $x \in W^J$ and $z \in W_J$. Moreover,
$\bar{G}=\bigsqcup\limits_{J \subset I} \bigsqcup\limits_{x \in
W^J, w \in W} [J, x, w]$. Springer proved the following result in
\cite{Spr1, 2.4}.

\proclaim{Theorem} Let $x \in W^J$, $x' \in W^K$, $w, w' \in W$.
Then $[K, x', w']$ is contained in $\overline{[J, x, w]}$ if and
only if $K \subset J$ and there exists $u \in W_K, v \in W_J \cap
W^K$ with $x v u \i \le x'$, $w' u \le w v$ and $l(w
v)=l(w)+l(v)$.
\endproclaim

As a consequence of the theorem, we have the following properties
which will be used later.

(1) For any $K \subset J$, $w \in W^J$ and $v \in W_J$, $[K, w v,
v] \subset \overline{[J, w, 1]}$.

(2) For any $J \subset I$, $w, w' \in W^J$ with $w \le w'$, then
$[J, w', 1] \subset \overline{[J, w, 1]}$.

\subhead 1.5 \endsubhead In this subsection, we recall some
results of \cite{Spr2}.

Let $\e$ be an indeterminate. Put $o=k[[\e]]$ and $K=k((\e))$. An
$o$-valued point of a $k$-variety $Z$ is a $k$-morphism $\g:
\Spec(o)@>>>Z$. We write $Z(o)$ for the set of all $o$-valued
points of $Z$. Similarly, we write $Z(K)$ for the set of all
$K$-valued points of $Z$. For $\g \in Z(o)$, we have that $\g(0)
\in Z$, where $0$ is the closed point of $\Spec(o)$.

By the valuative criterion of completeness (see \cite{EGA, Ch II,
7.3.8 \& 7.3.9}), for the complete $k$-variety $\bar{G}$, the
inclusion $o \hookrightarrow K$ induces a bijective from
$\bar{G}(o)$ onto $\bar{G}(K)$. Therefore, any $\g \in \bar{G}(K)$
defines a point $\g(0) \in \bar{G}$. In particular, any $\g \in
U(K)$ defines a point $\g(0) \in \bar{G}$. Here we regard $U(K)$
as a subset of $\bar{G}(K)$ in the natural way.

We have that $x \in \bar{U}$ if and only if there exists $\g \in
U(K)$ such that $\g(0)=x$ (see \cite{Spr2, 2.2}).

Let $Y$ be the cocharacter group of $T$. An element $\l \in Y$
defines a point in $T(k[\e, \e \i])$, hence a point $p_{\l}$ of
$T(K)$. Let $H \subset G(o)$ be the subgroup consisting of
elements $\g$ with $\g(0) \in B$. Then for $\g \in U(K)$, there
exists $\g_1, \g_2 \in H$, $w \in W$ and $\l \in Y$, such that
$\g=\g_1 \dot{w} p_{\l} \g_2$. Moreover, $w$ and $\l$ are uniquely
determined by $\g$ (see \cite{Spr2, 2.6}). In this case, we will
call $(w, \l)$ admissible. Springer showed that $(w, \l-w \i \l)$
is admissible for any dominant regular coweight $\l$ (see
\cite{Spr2, 3.1}).

For $\l \in Y$ and $x \in W$ with $x \i \cdot \l$ dominant, we
have that $p_{\l}(0)=(\dot{x}, \dot{x}) \cdot h_{I(x \i \l)}$,
where $I(x \i \l)$ is the set of simple roots orthogonal to $x \i
\l$ (see \cite{Spr2, 2.5}). If moreover, $(w, \l)$ is admissible,
then there exists some $t \in T$ such that $(U \times U) (\dot w
\dot x t, \dot x) \cdot h_{I(x \i \l)} \subset \bar{U}$.

\head 2. the partition of $Z_J$ \endhead

\subhead 2.1 \endsubhead We will follow the set-up of \cite{L4,
8.18}.

For any $J \subset I$, let $\cp^J$ be the set of parabolic
subgroups conjugate to $P_J$. We will write $\cb$ for
$\cp^{\varnothing}$. For $P \in \cp^J$, $Q \in \cp^K$ and $u \in
{^J W^K}$, we write $\po(P, Q)=u$ if there exists $g \in G$, such
that $^g P=P_J, {}^g Q={}^{\dot u} P_K$. For $J, J' \subset I$ and
$y \in {^{J'} W ^J}$ with $\Ad(y)J=J'$, define
$$\tz^y_J=\{(P, P', \g) \mid P \in \cp^J, P' \in \cp^{J'},
\g=U_{P'} g U_P, \po(P', {}^g P)=y\}$$ with the $G \times G$
action given by $(g_1, g_2) \cdot (P, Q, \g)=( {^{g_1}P},
{}^{g_2}P', g_2 \g g_1\i)$.

To $z=(P,P',\g)\in \tilde{Z^y_J}$, we associate a sequence
$(J_k,J'_k,u_k, y_k, P_k, P'_k, \g_k)_{k \ge 0}$ with $J_k,J'_k
\subset I$, $u_k \in W$, $y_k \in
{}^{J'_k}W^{J_k},\Ad(y_k)J_k=J'_k$, $P_k \in \cp_{J_k},P'_k \in
\cp_{J'_k},\g_k=U_{P'_k} g U_{P_k}$ for some $g \in G$ satisfies
$\po(P'_k, {}^g P_k)=u_k$. The sequence is defined as follows.
$$P_0=P,P'_0=P',\g_0=\g,J_0=J,J'_0=J',u_0=\po(P'_0,P_0),y_0=y.$$
Assume that $k \ge 1$, that $P_m,P'_m,\g_m,J_m,J'_m,u_m,y_m$ are
already defined for $m<k$ and that
$u_m=\po(P'_m,P_m),P_m\in\cp_{J_m},P'_m\in\cp_{J'_m}$ for $m<k$.
Let
$$J_k=J_{k-1} \cap \Ad(y_{k-1}\i u_{k-1})J_{k-1},
J'_k=J_{k-1} \cap \Ad(u_{k-1}\i y_{k-1})J_{k-1},$$
$$P_k=g_{k-1}\i({}^{g_{k-1}}P_{k-1})^{(P'_{k-1}{}^{P_{k-1}})}g_{k-1}\in\cp_{J_k},
P'_k=P_{k-1}^{P'_{k-1}}\in\cp_{J'_k}$$ where
$$g_{k-1} \in \g_{k-1} \hbox{ is such that } ^{g_{k-1}} P_{k-1} \hbox { contains some Levi of }
P_{k-1} \cap P'_{k-1},$$
$$u_k=\po(P'_k,P_k),y_k=u_{k-1}\i y_{k-1},\g_k=U_{P'_k}g_{k-1}U_{P_k}.$$

It is known that the sequence is well defined. Moreover, for
sufficient large $n$, we have that
$J_n=J'_n=J_{n+1}=J'_{n+1}=\cdots$ and $u_n=u_{n+1}=\cdots=1$. Now
we set $\b(z)=u_0 u_1 \cdots u_n$, $n \gg 0$. Then we have that
$\b(z) \in {}^{J'} W$. By \cite{L4, 8.18} and \cite{L3, 2.5}, the
sequence $(J_k, J'_k, u_k, y_k)_{k \ge 0}$ is uniquely determined
by $(J, \b(z), y)$.

The map $w \mapsto y w \i$ is a bijection between $W^J$ and $^{J'}
W$. For $w \in W^J$, set
$$\tz^y_{J, w}=\{z \in \tz^y_J \mid \b(z)=y w \i \}.$$

Then $(\tz^y_{J, w})_{w \in W^J}$ is a partition of $\tz^y_J$ into
locally closed $G$-stable subvarieties. For $w \in W^J$, let
$(J_k, J'_k, u_k, y_k)_{k \ge 0}$ be the sequence uniquely
determined by $(J, y w \i, y)$. Then $(P, P', \g)\mapsto(P_1,
P'_1, \g_1)$ define a $G$-equivariant map $\v: \tz^y_{J,
w}@>>>\tz^{y_1}_{J_1, u_0 \i w}$.

\subhead 2.2 \endsubhead Let $J \subset I$. Set $\tz_J=\tz^{w_0
w^J_0}_J$ and $J^*=\Ad(w_0 w^J_0) J$. For $w \in W^J$, set
$w_J=w_0 w^J_0 w \i$. The map $w \mapsto w_J$ is a bijection
between $W^J$ and $^{J^*} W$. For any $w \in W^J$, let
$$\tz_{J, w}=\{z \in \tz_J \mid \b(z)=w_J\}.$$

Then $(\tz_{J, w})_{w \in W^J}$ is a partition of $\tz_J$ into
locally closed $G$-stable subvarieties. Let $(J_k, J'_k, u_k,
y_k)_{k \ge 0}$ be the sequence determined by $(J, w_J, w_0
w^J_0)$ (see 2.1). Assume that $J_n=J'_n=J_{n+1}=J'_{n+1}=\cdots$
and $u_n=u_{n+1}=\cdots=1$. Set $v_0=w_J$ and $v_k=u_{k-1} \i
v_{k-1}$ for $k \in \bold N$. By \cite{L4, 8.18} and \cite{L3,
2.3}, we have $u_k \in {}^{J'_k} W^{J_k}$ and $u_{k+1} \in
W_{J_k}$ for all $k \ge 0$. Hence $v_{k+1} \in W_{J_k}$ for all $k
\ge 0$. Moreover, it is easy to see by induction on $k$ that
$y_k=v_k w$. In particular, $w=y_n \in {}^{J_n} W ^{J_n}$, $Ad(w)
J_n=J_n$ and $\dot w$ normalizes $B \cap L_{J_n}$. We have the
following result.

\proclaim{Lemma 2.3} Keep the notation of 2.2. Let $z=(P_J,
{}^{\dot w_J \i} P_{J^*}, \dot w_J \i U_{P_{J^*}} \dot w_J \dot w
b U_{P_J})$, where $b \in {}^{\dot w^{n-1} \dot v_n \i}
(U_{P_{J'_n}} \cap U_{J_{n-1}}) ^{\dot w^{n-2} \dot v_{n-1} \i}
(U_{P_{J'_{n-1}}} \cap U_{J_{n-2}}) \cdots ^{\dot v_1 \i}
(U_{P_{J'_1}} \cap U_{J_0})T$ or $b \in B$. Then $z \in \tz_{J,
w}$.
\endproclaim

Proof. For any $k$, set $P_k=P_{J_k}, P'_k={}^{\dot v_k \i}
P_{J'_k}$. Then

$$P_k \cap P'_k=P_{J_k} \cap {}^{\dot
v_{k+1} \i \dot u_k \i} P_{J'_k}={}^{\dot v_{k+1} \i} (P_{J_k}
\cap {}^{\dot u_k \i} P_{J'_k}).$$

Note that $u_k \i \in {}^{J_k} W ^{J'_k}$. Then $L_{J_k} \cap
{}^{\dot u_k \i} L_{J'_k}=L_{J_k \cap Ad(\dot u_k \i)
J'_k}=L_{J'_{k+1}}$. Thus $^{\dot v_{k+1} \i}
L_{J'_{k+1}}={}^{\dot v_{k+1} \i} (L_{J_k} \cap {}^{\dot u_k \i}
L_{J'_k})$ is a Levi factor of $P_k \cap P'_k$. Moreover, we have
$$\align P_k^{P'_k} &=P_{J_k}^{(^{\dot v_k \i} P_{J'_k})}={}^{\dot
v_{k+1} \i} (P_{J_k}^{(^{\dot u_k \i} P_{J'_k})})={}^{\dot v_{k+1}
\i} P_{J_k \cap Ad(\dot u_k \i) J'_k}={}^{\dot v_{k+1} \i}
P_{J'_{k+1}} \cr {P'_k}^{P_k} &={}^{\dot v_k \i}
(P_{J'_k}^{(^{\dot v_k} P_{J_k})})={}^{\dot v_k \i}
(P_{J'_k}^{(^{\dot u_k} P_{J_k})})={}^{\dot v_k \i} P_{J'_k \cap
Ad(\dot u_k) J_k} \cr &={}^{\dot v_k \i} P_{Ad(\dot y_k)(J_k \cap
Ad(\dot y_k \i \dot u_k) J_k)}={}^{\dot v_k \i} P_{Ad(\dot y_k)
J_{k+1}} \endalign$$

If $b \in B$, then set $g_k=\dot w b$, $\g_k=U_{P'_k} g_k U_{P_k}$
and $z_k=(P_k, P'_k, \g_k)$ for all $k$. In this case, $^{\dot
v_{k+1} \i} L_{J'_{k+1}}={}^{\dot w \dot y_{k+1} \i}
L_{J'_{k+1}}={}^{\dot w} L_{J_{k+1}} \subset {}^{\dot w}
P_k={}^{g_k} P_k$. Thus $^{g_k} P_k$ contains some Levi of $P_k
\cap P'_k$. Moreover, $$\eqalignno{g_k \i (^{g_k} P_k)^{(^{\dot
v_k \i} P_{Ad(\dot y_k) J_{k+1}})} g_k &=P_k^{(^{b \i \dot w \i
\dot v_k \i} P_{Ad(\dot y_k) J_{k+1}})}={}^{b \i} (P_k^{^{\dot y_k
\i} P_{Ad(\dot y_k) J_{k+1}}}) \cr &={}^{b \i} P_{J_k \cap Ad(\dot
y_k \i) Ad(\dot y_k) J_{k+1}}={}^{b \i}
P_{J_{k+1}}=P_{J_{k+1}}.}$$

Therefore, $\v(z_k)=z_{k+1}$.

If $b=(\dot w^{n-1} \dot v_n \i b_n \dot v_n \dot w^{-n+1}) \cdots
(\dot v_1 \i b_1 \dot v_1) (\dot w^n t \dot w^{-n})$, where $b_j
\in U_{P_{J'_j}} \cap U_{J_{j-1}}$ for $1 \le j \le n$ and $t \in
T$, then set
$$a_k=(\dot w^{n-k} \dot v_n \i b_n \dot
v_n \dot w^{-n+k}) \cdots (\dot v_k \i b_k \dot v_k) (\dot
w^{n+1-k} t \dot w^{-n-1+k}).$$

In this case, set $g_k=\dot w a_{k+1}$, $\g_k=U_{P'_k} g_k
U_{P_k}$ and $z_k=(P_k, P'_k, \g_k)$. For $j \ge 0$, $J_{j+1}=J_j
\cap \Ad(\dot y_{j+1} \i) J_j$ and $v_{j+1} \in W_{J_j}$. Thus
$^{\dot w} L_{J_{j+1}}={}^{\dot v_{j+1} \i \dot y_{j+1}}
L_{J_{j+1}} \subset {}^{\dot v_{j+1} \i} L_{J_j}=L_{J_j}$. Then
$^{\dot w^j \dot v_{k+j+1} \i} U_{J_{k+j}} \subset {}^{\dot w^j}
L_{J_{k+j}} \subset L_{J_k}$. So $a_{k+1} \in P_k$. Thus $^{g_k}
P_k={}^{\dot w} P_k$ contains some Levi of $P_{J_k} \bigcap ^{\dot
v_k \i} P_{J'_k}$. Moreover, $$g_k \i (^{g_k} P_k)^{(^{\dot v_k
\i} P_{Ad(\dot y_k) J_{k+1}})} g_k={}^{a_{k+1} \i} P_{J_{k+1}}.$$

Thus $\v(z_k)=(Q, Q', \g')$, where $Q={}^{a_{k+1} \i}
P_{J_{k+1}}$, $Q'={}^{\dot v_{k+1} \i} P_{J'_{k+1}}$ and
$\g'=U_{Q'} g_k U_Q$. Note that $^{\dot v_{k+1} \i}
U_{P_{J'_{k+1}}} \subset Q'$ and $T \subset Q'$. Moreover, for $j
\ge 1$, $^{\dot w^j \dot v_{k+j+1} \i} U_{J_{k+j}} \subset
{}^{\dot w^j} L_{J_{k+j}} \subset {}^{\dot w} L_{J_{k+1}}={}^{\dot
v_{k+1} \i \dot y_{k+1}} L_{J_{k+1}}={}^{\dot v_{k+1} \i}
L_{J'_{k+1}} \subset Q'$. Thus $a_{k+1} \in Q'$. Hence,
$z_{k+1}=(a_{k+1}, a_{k+1}) \cdot \v(z_k)$.

In both cases, $\v(z_k)$ is in the same $G$ orbit as $z_{k+1}$.
Thus $$\b(z)=\b(z_0)=u_1 \b(z_1)=\cdots=u_1 u_2 \cdots u_n=w_J.
\qed$$

\subhead Remark \endsubhead 1. From the proof of the case where $b
\in B$, we can see that $$\v^n(P_J, ^{\dot w_J \i} P_{J^*}, \dot
w_J \i U_{P_{J^*}} \dot w_J \dot w b U_{P_J})=(P_{J_n}, P_{J_n},
U_{P_{J_n}} \dot w b U_{P_{J_n}}).$$ This result will be used to
establish a relation between the $G$-stable pieces and the $B
\times B$-orbits.

2. The fact that $(P_J, {}^{\dot w_J \i} P_{J^*}, \dot w_J \i
U_{P_{J^*}} \dot w_J \dot w b U_{P_J})$ is contained in $\tz_{J,
w}$ for any $b \in {}^{\dot w^{n-1} \dot v_n \i} (U_{P_{J'_n}}
\cap U_{J_{n-1}}) ^{\dot w^{n-2} \dot v_{n-1} \i}
(U_{P_{J'_{n-1}}} \cap U_{J_{n-2}}) \cdots ^{\dot v_1 \i}
(U_{P_{J'_1}} \cap U_{J_0})T$ plays an important role in section
3. We will discuss about it in more detail in 3.1.

\subhead 2.4 \endsubhead  Let $(J_n, J'_n, u_n, y_n)_{n \ge 0}$ be
the sequence that is determined by $w_J$ and $w_0 w^J_0$. Assume
that $J_n=J'_n=J_{n+1}=J'_{n+1}=\cdots$ and
$u_n=u_{n+1}=\cdots=1$. Then $z \mapsto \v^n(z)$ is a
$G$-equivariant morphism from $\tz_{J, w}$ to $\tz^w_{J_n, 1}$ and
induces a bijection from the set of $G$-orbits on $\tz_{J, w}$ to
the set of $G$-orbits on $\tz^w_{J_n, 1}$.

Set $\tl_{J, w}=L_{J_n}$ and $\tc_{J, w}=\dot w \tl_{J, w}$. Let
$N_G(\tl_{J, w})$ be the normalizer of $\tl_{J, w}$ in $G$. Then
$\tc_{J, w}$ is a connected component of $N_G(\tl_{J, w})$ and
$\tz^w_{J_n, 1}$ is a fibre bundle over $\cp^{J_n}$ with fibres
isomorphic to $\tc_{J, w}$. There is a natural bijection between
$\tc_{J, w}$ and $F=\{z=(P_{J_n}, P_{J_n}, \g_n) \mid z \in
\tz^w_{J_n, 1}\}$ under which the action of $\tl_{J, w}$ on
$\tc_{J, w}$ by conjugation corresponds to the action of
$P_{J_n}/U_{P_{J_n}}$ on $F$ by conjugation. Therefore, we obtain
a canonical bijection the set of $G$-stable subvarieties of
$\tz_{J, w}$ and the set of $\tl_{J, w}$-stable subvarieties of
$\tc_{J, w}$ (see \cite{L4, 8.21}). Moreover, a $G$-stable
subvariety of $\tz_{J, w}$ is closed if and only if the
corresponding $\tl_{J, w}$-stable subvariety of $\tc_{J, w}$ is
closed. By the remark 1 of 2.3, for any $b \in B \cap \tl_{J, w}$,
the $G$-orbit that contains $(P_J, {}^{\dot w_J \i} P_{J^*}, \dot
w b)$ corresponds to the $\tl_{J, w}$-orbit that contains $\dot w
b$ via the bijection.

\subhead 2.5 \endsubhead Since $G$ is adjoint, the center of
$P/U_P$ is connected for any parabolic subgroup $P$. Let $H_P$ be
the inverse image of the (connected) center of $P/U_P$ under
$P@>>>P/U_P$. We can regard $H_P/U_P$ as a single torus $\D_J$
independent of $P$. Now $\D_J$ acts (freely) on $\tz_J$ by $\d:
(P,P',\g) \mapsto (P,P',\g z)$ where $z\in H_P$ represents
$\d\in\D_J$. The action of $G$ on $\tz_J$ commutes with the action
of $\D_J$ and induces an action of $G$ on $\D_J \setminus \tz_J$.
There exists a $G$-equivariant isomorphism from $Z_J$ to $\D_J
\setminus \tz_J$ which sends $(g_1, g_2) \cdot h_J$ to $(^{g_2}
P_J, {}^{g_1} P^-_J, U_{^{g_1} P^-_J} g_1 g_2 \i H_{^{g_2} P_J})$.
We will identify $Z_J$ with $\D_J \setminus \tz_J$.

It is easy to see that $\D_J(\tz_{J, w})=\tz_{J, w}$. Set $Z_{J,
w}=\D_J \setminus \tz_{J, w}$. Then $$Z_J=\bigsqcup_{w \in W^J}
Z_{J, w}.$$

Moreover, we may identify $\D_J$ with a closed subgroup of the
center of $\tl_{J, w}$. Set $L_{J, w}=\tl_{J, w}/ \D_J$ and $C_{J,
w}=\tc_{J, w} / \D_J$. Thus we obtain a bijection between the set
of $G$-stable subvarieties of $Z_{J, w}$ and the set of $L_{J,
w}$-stable subvarieties of $C_{J, w}$ (see \cite{L4, 11.19}).
Moreover, a $G$-stable subvariety of $Z_{J, w}$ is closed if and
only if the corresponding $L_{J, w}$-stable subvariety of $C_{J,
w}$ is closed and for any $b \in B \cap \tl_{J, w}$, the $G$-orbit
that contains $(P_J, {}^{\dot w_J \i} P_{J^*}, \dot w b)$
corresponds to the $L_{J, w}$-orbit that contains $\dot w b \D_J$
via the bijection.

\proclaim{Proposition 2.6} For any $w \in W^J$, $Z_{J, w}=G_{diag}
\cdot [J, w, 1]$.
\endproclaim

Proof. By 2.3, $(\dot w, b) \cdot h_J \in Z_{J, w}$ for all $b \in
B$. Since $Z_{J, w}$ is $G$-stable, $G_{diag} [J, w. 1] \subset
Z_{J, w}$.

For any $z \in Z_{J, w}$, let $C$ be the $L_{J, w}$-stable
subvariety corresponding to $G_{diag} \cdot z$ and let $c$ be an
element in $\tc_{J, w}$ such that $c \D_J \in C$. By 2.2, $\dot w$
normalizes $B \cap \tl_{J, w}$. Thus $c$ is $\tl_{J, w}$-conjugate
to an element of $\dot w (B \cap \tl_{J, w})$. Therefore, $z$ is
$G$-conjugate to $(\dot w, b) \cdot h_J$ for some $b \in B \cap
\tl_{J, w}$. The proposition is proved. \qed

\proclaim{Proposition 2.7} For any $w \in W^J$, $\overline{Z_{J,
w}}=\overline{G_{diag} (\dot w T, 1) \cdot h_J}$.
\endproclaim

Proof. Since $(\dot w T, 1) \cdot h_J \subset Z_{J, w}$ and
$\overline{Z_{J, w}}$ is a $G$-stable closed variety, we have that
$\overline{G_{diag} (\dot w T, 1) \cdot h_J} \subset
\overline{Z_{J, w}}$.

Set $X=\{(\dot w t, u) \cdot h_J \mid t \in T, u \in U\}$. For any
$u \in {}^{\dot w} U_J$ and $t \in T$, we have that $Ad(\dot w t)
\i u \in U_J$ and $u \in {}^{\dot w} U_J \subset U$. Consider the
map $\phi: ^{\dot w} U_J \times T@>>>X$ defined by $\phi(u, t)=(u,
u) (\dot w t, 1) \cdot h_J=(\dot w t, (\dot w t) \i u \dot w t u
\i) \cdot h_J$, for $u \in {}^{\dot w} U_J, t \in T$.

It is easy to see that there is an open subset $T'$ of $T$, such
that the restriction of $\phi$ to $^{\dot w} U_J \times T'$ is
injective. Note that
$\dim(X)=\dim(T)+\dim(U/U_{P_J})=\dim(T)+\dim(U_J)=\dim(^{\dot w}
U_J \times T)$. Then the image of $\phi$ is dense in $X$. The
proposition is proved. \qed

\subhead Remark \endsubhead This argument was suggested by the
referee.

\subhead 2.8 \endsubhead For $w \in W$, denote by $\supp(w)$ the
set of simple roots whose associated simple reflections occur in a
reduced expression of $w$. An element $w \in W$ is called a
Coxeter element if it is a product of the simple reflections, in
some order, or in other words, $|\supp(w)|=l(w)=|I|$. We have the
following properties.

\proclaim{Proposition 2.9} Fix $i \in I$. Then all the Coxeter
elements are conjugate under elements of $W_{I-\{i\}}$.
\endproclaim

Proof. Let $c, c'$ be Coxeter elements. We say that $c'$ can be
obtained from $c$ via a cyclic shift if $c=s_{i_1} s_{i_2} \cdots
s_{i_n}$ is a reduced expression and $c'=s_{i_1} c s_{i_1}$. It is
known that for any Coxeter elements $c, c'$, there exists a finite
sequences of Coxeter elements $c=c_0, c_1, \ldots, c_m=c'$ such
that $c_{k+1}$ can be obtained from $c_k$ via a cyclic shift (see
\cite{Bo, p. 116, Prop. 1}).

Now assume that $c=s_{i_1} s_{i_2} \cdots s_{i_n}$ is a reduced
expression of a Coxeter element. If $i_1 \neq i$, then $s_{i_1} c
s_{i_1}$ and $c$ are conjugated by $s_{i_1} \in W_{I-\{i\}}$. If
$i_1=i$, then $s_{i_1} c s_{i_1}=s_{i_2} s_{i_3} \cdots s_{i_n} c
(s_{i_2} s_{i_3} \cdots s_{i_n}) \i$. Therefore, if a Coxeter
element can be obtained from another Coxeter element via a cyclic
shift, then they are conjugated by elements of $W_{I-\{i\}}$. The
proposition is proved. \qed

\subhead Remark \endsubhead The proof of \cite{loc. cit} also can
be used to prove this proposition.

\proclaim{Proposition 2.10} Let $J \subset I$ and $w \in W^J$ with
$\supp(w)=I$. Then there exist a Coxeter element $w'$, such that
$w' \in W^J$ and $w' \le w$.
\endproclaim

Proof. We prove the statement by induction on $l(w)$.

Let $i \in I$ with $s_i w<w$. Then $s_i w \in W^J$. If $\supp(s_i
w)=I$, then the statement holds by induction hypothesis on $s_i
w$. Now assume that $\supp(s_i w)=I-\{i\}$. By induction, there
exists a Coxeter element $w'$ of $W_{I-\{i\}}$, such that $w' \in
W^{J-\{i\}}$ and $w' \le s_i w$. Then $s_i w'$ is a Coxeter
element of $w$ and $s_i w' \le w$.

Since $w' \in W_{I-\{i\}}$, $(w') \i \a_i$ is either $\a_i$ or a
non-simple positive root. We also have that $w'$ is a Coxeter
element of $W_{I-\{i\}}$. Thus if $(w') \i \a_i=\a_i$, then
$<\a_i, \a^\vee_j>=0$ for all $j \neq i$. It contradicts the
assumption that $G$ is simple. Hence $(w') \i \a_i$ is a
non-simple positive root. Note that if $s_i w' \notin W^J$, then
$s_i w'=w' s_j$ for some $j \in J$, that is, $(w') \i \a_i=\a_j$.
Therefore, $s_i w' \in W^J$. The proposition is proved. \qed

\proclaim{Corollary 2.11} Let $i \in I$, $J=I-\{i\}$ and $w$ be a
Coxeter element of $W$ with $w \in W^J$. Then $\bigsqcup_{K
\subset J} \bigsqcup_{w' \in W^K, \supp(w')=I} Z_{K, w'} \subset
\overline{Z_{J, w}}$.
\endproclaim

Proof. By 1.4, $[K, w v, v] \subset \overline{[J, w, 1]}$ for $K
\subset J$ and $v \in W_J$. Since $\overline{Z_{J, w}}$ is
$G$-stable, $(\dot v \i \dot w \dot v T, 1) \cdot h_K \subset
\overline{Z_{J, w}}$. By 2.9, $(\dot w' T, 1) \cdot h_K \subset
\overline{Z_{J, w}}$ for all Coxeter element $w'$. By 2.7, $Z_{K,
w'} \subset \overline{Z_{J, w}}$ for all Coxeter element $w'$ with
$w' \in W^K$. For any $u \in W^K$ with $\supp(u)=I$, there exists
a Coxeter element $w'$, such that $w' \in W^K$ and $w' \le u$.
Thus by 1.4, we have that $[K, u, 1] \subset \overline{Z_{J, w}}$.
By 2.6, $Z_{K, u} \subset \overline{Z_{J, w}}$. The corollary is
proved. \qed

\subhead Remark \endsubhead In 4.4, we will show that the equality
holds.

\head 3. Some combinatorial results \endhead

\subhead 3.1 \endsubhead Fix $i \in I$. Define subsets $I_k$ of
$I$ for all $k \in \bold N$ in the following way. Set $I_1=\{i\}$.
Assume that $I_k$ is already defined. Set
$$I_{k+1}=\{\a_j \mid j \in I-\cup_{l=1}^k I_l, <\a^\vee_j, \a_m>
\neq 0 \hbox{ for some } m \in I_k\}.$$

It is easy to see that if $j_1, j_2 \in I_k$ with $j_1 \neq j_2$,
then $<\a_{j_1}, \a^\vee_{j_2}>=0$. Thus $s_{I_k}=\prod_{j \in
I_k} s_j$ is well-defined. For sufficiently large $n$, we have
$I_n=I_{n+1}=\cdots=\varnothing$ and
$s_{I_n}=s_{I_{n+1}}=\cdots=1$. Now set $w_k=s_{I_n} s_{I_{n-1}}
\cdots s_{I_k}$ for $k \in \bold N$. We will write $w^J$ for
$w_1$. Set $J_{-1}=I$ and $J_0=J=I-\{i\}$. Then $w^J$ is a Coxeter
element and $w^J \in W^J$. Let $(J_n, J'_n, u_n, y_n)$ be the
sequence determined by $w^J$ and $w_0 w^J_0$. Then we can show by
induction that for $k \ge 0$, $J_k=J_{k-1}-I_{k+1}$,
$u_k=w^{J_{k-1}}_0 w^{J_k}_0 s_{I_{k+1}} w^{J_{k+1}}_0 w^{J_k}_0$,
$y_k=w^{J_{k-1}}_0 w^{J_k}_0 s_{I_k} s_{I_{k-1}} \cdots s_{I_1}$
and $J'_k=Ad(y_k) J_k$. In particular, $J_n=\varnothing$. Thus
$\tl_{J, w^J}=T$ and $\tc_{J, w^J}=\dot w^J T$. Since $w$ is a
Coxeter element, the homomorphism $T@>>>T$ sending $t \in T$ to
$(\dot w^J)^{-1} t \dot w^J t \i$ is surjective. Thus $\tl_{J,
w^J}$ acts transitively on $\tc_{J, w^J}$. By 2.5, $G$ acts
transitively on $Z_{J, w^J}$.

For $k \in \bold N$, we set $v_k=w^{J_{k-1}}_0 w^{J_k}_0 w_{k+1}
\i$. Then it is easy to see that $$^{\dot v_k \i} (U_{P_{J'_k}}
\cap U_{J_{k-1}})={}^{w_{k+1}}(U_{P^-_{J_k}} \cap
U^-_{J_{k-1}}).$$ Therefore by 2.3, for $b \in {}^{w^{n-1}
w_{n+1}} (U_{P^-_{J_n}} \cap U^-_{J_{n-1}}) \cdots ^{w_2}
(U_{P^-_{J_1}} \cap U^-_{J_0}) T$, we have that $(\dot w^J b, 1)
\cdot h_J \in Z_{J, w^J}$.

In the rest of this section, we will keep the notations of $J$,
$J_k$, $w^J$ and $w_k$ as above. We will prove the following
statement.

\proclaim{Proposition} Let $X$ be a closed subvariety of $\bar{G}$
satisfying the following condition: for any admissible pair $(w,
\l)$ and $x \in W$ with $x \i \l$ is dominant, there exist some $t
\in T$, such that $G_{diag} (U \times U) (\dot w \dot x t, \dot x)
\cdot h_{I(x \i \l)} \subset X$. Then $Z_{J, w^J} \subset X$.
\endproclaim

An example of such $X$ is $\bar{\cu}$. There are some other
interesting examples, which we will discuss in 4.5. The proof is
based on case-by-case checking.


\subhead Remark \endsubhead The outline of the case-by-case
checking is as follows.

For $\l \in Y$, we write $\l \ge 0$ if $\l \in \sum_{l \in I}
\bold R_{\ge 0} \a^\vee_l$.

We start with the fundamental coweight $\o^\vee_i$. Find $x \in W$
that satisfies the conditions (1) $x \o^\vee_i \ge 0$ and (2) for
$l \in I$, either $(s_l-1) x \o^\vee_i \ge 0$ or $s_l x \o^\vee_i
\ngeqslant 0$. Such $x$ always exists, as we will see by
case-by-case checking. The elements $x \o^\vee_i$ that we obtain
in this way are not unique, in general. Fortunately, there always
exists some $x \in W$ that satisfies the conditions (1) and (2)
and allows us to do the procedures that we will discuss below.

In the rest of the remark, we fix such $x$. Since $x \o^\vee_i \in
Y$, there exists $n \in \bold N$, such that $n x \o^\vee_i$ is
contained in the coroot lattice. Set $\l=n x \o^\vee_i$. Now we
can find $v \in W$ such that $(v, \l)$ is admissible. (In
practice, we find $v \in W$ with $l(v)=|\supp(v)|$ and $-v \l \ge
0$. Then we can use lemma 3.2 to check that if $(v, \l)$ is
admissible.) By the assumption on $X$, $G_{diag} (U \times U)
(\dot v \dot x t, \dot x) \cdot h_J \subset X$ for some $t \in T$.

In some cases, $x \i v x=w_J$. Since $w_J$ is a Coxeter element,
$(\dot w_J T, 1) \cdot h_J=T_{diag} (\dot w_J t, 1) \cdot h_J
\subset X$. By 2.7, $Z_{J, w_J} \subset X$.

In other cases, the situation is more complicated. We need to
choose some $u \in U$, such that $(u \dot v \dot x t, \dot x)
\cdot h_J \in Z_{J, w_J}$. This is the most difficult part of the
case-by-case checking. The lemma 3.3 and lemma 2.3 will be used to
overcome the difficulties.

${}$

Throughout this section, we will use the same labelling of Dynkin
diagram as in \cite{Bo}. For $a, b \in I$, we denote by $s_{[a,
b]}$ the element $s_b s_{b-1} \cdots s_a$ of the Weyl group $W$
and $\dot s_{[a, b]}=\dot s_b \dot s_{b-1} \cdots \dot s_a$. (If
$b<a$, then $s_{[a, b]}=1$ and $\dot s_{[a, b]}=1$.)

\proclaim{Lemma 3.2} Let $x=s_{i_1} s_{i_2} \cdots s_{i_n}$ with
$|\supp(x)|=n$. Then $(1-x \i) \o^\vee_k=0$ if $k \notin \{i_1,
i_2, \ldots, i_n\}$ and $(1-x \i) \o^\vee_{i_j}=s_{i_n}
s_{i_{n-1}} \cdots s_{i_{j+1}} \a^\vee_{i_j}$. Thus $(x, \l)$ is
admissible for all $\l \in \sum_{j=1}^n \bold N s_{i_n}
s_{i_{n-1}} \cdots s_{i_{j+1}} \a^\vee_{i_j}$.
\endproclaim

The lemma is a direct consequence of \cite{Bo, p. 226, Ex. 22a},
which was pointed out to me by the referee.

\proclaim{Lemma 3.3} Let $w, x, y_1, y_2 \in W$ and $t \in T$.
Assume that $y_1=s_{i_1} s_{i_2} \cdots s_{i_l}$, $y_2=s_{i_{l+1}}
s_{i_{l+2}} \cdots s_{i_{l+k}}$ with $k+l=|\supp(y_1 y_2)|$. If
moreover, $<\a^\vee_{i_{l_1}}, \a_{i_{l_2}}>=0$ for all $1 \le
l_1<l_2 \le l$ and $(1-y_1 y_2) x \o^\vee_i, (1-y_1) w \o^\vee_i
\in \sum_{j=1}^k \bold R_{>0} \a^\vee_{i_j}$, then there exists $u
\in U_{-w \i \a_{i_{l+1}}} U_{-w \i \a_{i_{l+2}}} \cdots U_{-w \i
\a_{i_{l+k}}}$ such that $(\dot x \i \dot w u t, 1) \cdot h_J \in
G_{diag} (U \times U) (\dot w t, \dot y_1 \dot y_2 \dot x) \cdot
h_J$.
\endproclaim

Proof. We have that $(1-y_1 y_2) x
\o^\vee_i=\sum_{j=1}^{k+l}(1-s_{i_j}) s_{i_{j+1}} \cdots
s_{i_{l+k}} x \o^\vee_i$. Note that $i_1, i_2, \ldots, i_{k+l}$
are distinct and $(1-s_{i_j}) s_{i_{j+1}} \cdots s_{i_{l+k}} x
\o^\vee_i \in \bold R \a^\vee_{i_j}$ for all $j$. Hence
$(1-s_{i_j}) s_{i_{j+1}} \cdots s_{i_{l+k}} x \o^\vee_i \in \bold
R_{>0} \a^\vee_{i_j}$ for all $j$, i. e., $<s_{i_{j+1}} \cdots
s_{i_k} x \o^\vee_i, \a_{i_j}> \in \bold R_{>0}$. Therefore $\dot
x \i \dot s_{i_{l+k}} \i \cdots \dot s_{i_{j+1}} \i U_{\a_{i_j}}
\dot s_{i_{j+1}} \cdots \dot s_{i_{l+k}} \dot x \subset U_{P_J}$.
Similarly, we have that $\dot w \i U_{-\a_{i_j}} \dot w \in
U_{P^-_J}$ for $j \le l$.

There exists $u_j \in U_{\a_{i_j}}$ and $u'_j \in U_{-\a_{i_j}}$
such that $u_j \dot s_{i_j} u_j=u'_j$. Note that $u'_1 u'_2 \cdots
u'_{l+k-1} \in L_{I-\{i_{l+k}\}}$, $u_{l+k} \in
U_{P_{I-\{i_{l+k}\}}}$ and $\dot x \i u_{l+k} \dot x \subset
U_{P_J}$. Thus
$$\eqalignno{u'_1 u'_2 \cdots u'_{l+k} \dot
x &=u'_1 u'_2 \cdots u'_{l+k-1} u_{l+k} \dot s_{i_k} u_{l+k} \dot
x \in U_{P_{I-\{i_k\}}} u'_1 u'_2 \cdots u'_{l+k-1} \dot s_{i_k}
\dot x U_{P_J} \cr & \subset U u'_1 u'_2 \cdots u'_{l+k-1} \dot
s_{i_k} \dot x U_{P_J}.}$$

We can show in the same way that $u'_1 u'_2 \cdots u'_{l+k} \dot x
\in U \dot y_1 \dot y_2 \dot x U_{P_J}$. Therefore, $(\dot w t,
u'_1 u'_2 \cdots u'_{l+k} \dot x) \cdot h_J \in (U \times U) (\dot
w t, \dot y_1 \dot y_2 \dot x) \cdot h_J$. Set $u=\dot w \i
u'_{l+1} u'_{l+2} \cdots u'_{l+k} \dot w$ and $u'=t \i \dot w \i
(u'_1 u'_2 \cdots u'_l) \i \dot w t \in U_{P^-_J}$. Then
$$\eqalignno{(\dot x \i \dot w u t, 1) \cdot h_J &=(\dot x \i \dot w u t u',
1) \cdot h_J=\bigl(\dot x \i (u'_1 u'_2 \cdots u'_{l+k}) \i \dot w
t, 1) \cdot h_J \cr & \in G_{diag} (U \times U) (\dot w t, \dot
y_1 \dot y_2 \dot x) \cdot h_J. \qed}$$

\subhead 3.4 \endsubhead In subsection 3.4 to subsection 3.7, we
assume that $G$ is $PGL_n(k)$. Without loss of generality, we
assume that $i \le n/2$. In this case, $w^J=s_{[i+1, n-1]} s_{[1,
i]} \i$. For any $a \in \bold R$, we denote by $[a]$ the maximal
integer that is less than or equal to $a$.

For $1 \le j \le i$, set $a_j=[(j-1)n/i]$. For convenience, we
will set $a_{i+1}=n-1$. Note that for $j \le i-1$, $a_{j+1}-a_j=[j
n/i]-[(j-1)n/i] \ge [n/i] \ge 2$. Therefore, we have that
$0=a_1<a_1+1<a_2<a_2+1< \cdots <a_i<a_i+1 \le a_{i+1}=n-1$. Now
set $b_0=0$. For $k \in \{1, 2, \ldots, n-1\}-\{a_2, a_3, \ldots,
a_i\}-\{a_2+1, a_3+1, \ldots, a_i+1\}$, set $b_k=i$. For $j \in
\{2,3, \ldots, i\}$, set $b_{a_j}=(j-1)n-i a_k$ and
$b_{a_j+1}=i-b_{a_j}$. In particular, $b_{n-1}=i$.

Now set $v=s_{[a_1+1, a_2-\d_{b_{a_2}, 0}]} s_{[a_2+1,
a_3-\d_{b_{a_3}, 0}]} \cdots s_{[a_i+1, a_{i+1}-\d_{b_{a_{i+1}},
0}]}$, where $\d_{a, b}$ is the Kronecker delta. Set
$v_j=s_{[a_j+1, a_{j+1}]} s_{[a_{j+1}+1, a_{j+2}]} \cdots
s_{[a_i+1, a_{i+1}]}$ for $1 \le j \le i$. Set $\l=\sum_{j=1}^i
\sum_{k=1}^{a_{j+1}-a_j} b_{a_j+k} (s_{[a_j+1, a_j+k-1]} v_{j+1})
\i \a^\vee_{a_j+k}$. It is easy to see that for $1 \le a \le b \le
n-1$ and $1 \le k \le n-1$,

$$s_{[a, b]} \a^\vee_k=\cases \sum_{l=a-1}^b \a^\vee_l, & \qquad \hbox{ if
} k=a-1; \cr -\sum_{l=a}^b \a^\vee_l, & \qquad \hbox{ if } k=a;
\cr \a^\vee_{k-1}, & \qquad \hbox{ if } a<k \le b; \cr
\a^\vee_b+\a^\vee_{b+1}, & \qquad \hbox{ if } k=b+1; \cr
\a^\vee_k, & \qquad \hbox{ otherwise }. \cr
\endcases$$

If $b_{a_j+k} \neq 0$, then $(s_{[a_j+1, a_j+k-1]} s_{[a_{j+1}+1,
a_{j+2}-\d_{b_{a_{j+2}}, 0}]} \cdots s_{[{a_i}+1, a_{i+1}]}) \i
\a^\vee_{a_j+k}=(s_{[a_j+1, a_j+k-1]} v_{j+1}) \i
\a^\vee_{a_j+k}$. By 3.2, $(v, \l)$ is admissible.

We have that $$\eqalignno{\l &=\sum_{j=1}^i
\sum_{k=1}^{a_{j+1}-a_j-1} b_{a_j+k} v_{j+1} \i s_{[a_j+1,
a_j+k-1]} \i \a^\vee_{a_j+k}+\sum_{j=1}^i b_{a_{j+1}} v_{j+1} \i
s_{[a_j+1, a_{j+1}-1]} \i \a^\vee_{a_{j+1}} \cr &=\sum_{j=1}^i
\sum_{k=1}^{a_{j+1}-a_j-1} \sum_{l=1}^k b_{a_j+k}
\a^\vee_{a_j+l}+\sum_{j=1}^{i-1} b_{a_{j+1}}
\sum_{l=1}^{a_{j+1}-a_j+1} \a^\vee_{a_j+l}+b_{a_{i+1}}
\sum_{l=1}^{a_{i+1}-a_i} \a^\vee_{a_i+l} \cr &=\sum_{j=1}^i
\sum_{k=1}^{a_{j+1}-a_j} \sum_{l=1}^k b_{a_j+k}
\a^\vee_{a_j+l}+\sum_{j=1}^{i-1} b_{a_{j+1}} \a^\vee_{a_{j+1}+1}
\cr &=\sum_{j=1}^i \sum_{l=1}^{a_{j+1}-a_j}
\sum_{k=l}^{a_{j+1}-a_j} b_{a_j+k}
\a^\vee_{a_j+l}+\sum_{j=1}^{i-1} b_{a_{j+1}} \a^\vee_{a_{j+1}+1}
\cr &=\sum_{j=1}^i \sum_{l=2}^{a_{j+1}-a_j} \bigl(
(a_{j+1}-a_j-l)i+b_{a_{j+1}} \bigr)
\a^\vee_{a_j+l}+\bigl((a_2-1)i+b_{a_2} \bigr) \a^\vee_1 \cr &
\qquad + \sum_{j=2}^i \bigl(b_{a_j}+(a_{j+1}-a_j-2) i
+b_{a_{j+1}}+b_{a_j+1} \bigr) \a^\vee_{a_j+1} \cr &=\sum_{j=1}^i
\sum_{l=1}^{a_{j+1}-a_j} \bigl( (a_{j+1}-a_j-l)i+b_{a_{j+1}}
\bigr) \a^\vee_{a_j+l}=n x \o^\vee_i.}$$

Note that $a_j \ge j$ for $j \ge 2$. Set $x_i=1$ and $x_j=s_{[j+1,
a_{j+1}]} s_{[j+2, a_{j+2}]} \cdots s_{[i, a_i]}$ for $1 \le j \le
i-1$. If $j=1$, we will simply write $x$ for $x_1$.

\proclaim{Lemma 3.5} For $1 \le j \le i$, we have that
$$\eqalignno{n x_j \o^\vee_i &=\sum_{l=1}^{j-1} l(n-i) \a^\vee_l+\sum_{l=j}^{a_{j+1}} \bigl(
j n-i l \bigr) \a^\vee_l \cr & \qquad +\sum_{k=j+1}^i
\sum_{l=1}^{a_{k+1}-a_k} \bigl( (a_{k+1}-a_k-l) i+b_{a_{k+1}}
\bigr) \a^\vee_{a_k+l}.}$$

In particular, $n x \o^\vee_i=\sum_{j=1}^i
\sum_{l=1}^{a_{j+1}-a_j} \bigl( (a_{j+1}-a_j-l)i+b_{a_{j+1}}
\bigr) \a^\vee_{a_j+l}$.
\endproclaim

Proof. We argue by  induction on $j$. Note that $n
\o^\vee_i=\sum_{l=1}^{i-1} l(n-i) \a^\vee_l+\sum_{l=i}^{n-1}
i(n-l) \a^\vee_l$. Thus the lemma holds for $j=i$.

Note that $j n-i (a_j+l)=j n -i
a_{j+1}+i(a_{j+1}-a_j-l)=b_{a_{j+1}}+i(a_{j+1}-a_j-l)$. Assume
that the lemma holds for $j$. Then
$$\eqalignno{& n x_{j-1} \o^\vee_i=s_{[j, a_j]}
\sum_{l=1}^{j-1} l(n-i) \a^\vee_l+s_{[j, a_j]}
\sum_{l=j}^{a_{j+1}} (j n-i l) \a^\vee_1\cr & \qquad \qquad \qquad
\qquad +s_{[j, a_j]} \sum_{k=j+1}^i \sum_{l=1}^{a_{k+1}-a_k}
\bigl( (a_{k+1}-a_k-l) i+b_{a_{k+1}} \bigr) \a^\vee_{a_k+l} \cr
&=\sum_{l=1}^{j-2} l(n-i) \a^\vee_l+(j-1)(n-i) \sum_{l=j-1}^{a_j}
\a^\vee_l-j(n-i) \sum_{l=j}^{a_j} \a^\vee_l+\sum_{l=j+1}^{a_j} (j
n-i l) \a^\vee_{l-1} \cr & \qquad \qquad +\bigl( j n-i (a_j+1)
\bigr) (\a^\vee_{a_j}+\a^\vee_{a_{j+1}})+\sum_{l=a_j+2}^{a_{j+1}}
(j n-i l) \a^\vee_l \cr & \qquad \qquad +\sum_{k=j+1}^i
\sum_{l=1}^{a_{k+1}-a_k} \bigl( (a_{k+1}-a_k-l) i+b_{a_{k+1}}
\bigr) \a^\vee_{a_k+l} \cr &=\sum_{l=1}^{j-2} l(n-i)
\a^\vee_l+(j-1)(n-i) \sum_{l=j-1}^{a_j} \a^\vee_l-j(n-i)
\sum_{l=j}^{a_j} \a^\vee_l+\sum_{l=j+1}^{a_j} (j n-i l)
\a^\vee_{l-1} \cr & \qquad \qquad +\bigl( j n-i (a_j+1) \bigr)
\a^\vee_{a_j}+\sum_{k=j}^i \sum_{l=1}^{a_{k+1}-a_k} \bigl(
(a_{k+1}-a_k-l) i+b_{a_{k+1}} \bigr) \a^\vee_{a_k+l} \cr
&=\sum_{l=1}^{j-2} l(n-i) \a^\vee_l+(j-1)(n-i)
\a^\vee_{j-1}+\sum_{l=j}^{a_j} \bigl( (j-1)(n-i)-j(n-i)+j n-i(l+1)
\bigr) \a^\vee_l\cr & \qquad \qquad +\sum_{k=j}^i
\sum_{l=1}^{a_{k+1}-a_k} \bigl( (a_{k+1}-a_k-l) i+b_{a_{k+1}}
\bigr) \a^\vee_{a_k+l} \cr &=\sum_{l=1}^{j-2} l(n-i)
\a^\vee_l+\sum_{l=j-1}^{a_j} \bigl( (j-1)n-i l \bigr) \a^\vee_l
+\sum_{k=j}^i \sum_{l=1}^{a_{k+1}-a_k} \bigl( (a_{k+1}-a_k-l)
i+b_{a_{k+1}} \bigr) \a^\vee_{a_k+l}.}$$

Thus the lemma holds for $j$. \qed

\proclaim{Lemma 3.6} We have that $x \i v_1 x=w^J$.
\endproclaim

Proof. If $a_j \ge j+1$, then $s_{[j+1, a_{j+1}]} \i s_{[a_j+1,
a_{j+1}]}=s_{[j+1, a_j]} \i$. If $j \ge 2$ and $a_j<j+1$, then
$j=2$, $a_j=2$ and $s_{[3, a_3]} \i s_{[a_2+1, a_3]}=1=s_{[3,
a_2]} \i$. In conclusion, $s_{[j+1, a_{j+1}]} \i s_{[a_j+1,
a_{j+1}]}=s_{[j+1, a_j]} \i$ for $j \ge 2$. Moreover, $s_{[2,
a_2]} \i s_{[a_1+1, a_2]}=s_1$. Thus

$$
s_{[2, a_2]} \i v_1 s_{[2, a_2]}=s_{[2, a_2]} \i s_{[a_1+1, a_2]}
v_2 s_{[2, a_2]}=s_1 v_2 s_{[2, a_2]}=v_2 s_1 s_{[2, a_2]}=v_2
s_{[3, a_2]} s_1 s_2.$$
$$\eqalignno{s_{[j+1, a_{j+1}]} \i v_j
s_{[j+1, a_j]} s_{[1, j]} \i s_{[j+1, a_{j+1}]} &=s_{[j+1,
a_{j+1}]} \i s_{[a_j+1, a_{j+1}]} v_{j+1} s_{[j+1, a_j]} s_{[1,
j]} \i s_{[j+1, a_{j+1}]} \cr &=s_{[j+1, a_j]} \i v_{j+1} s_{[j+1,
a_j]} s_{[1, j]} \i s_{[j+1, a_{j+1}]} \cr &=v_{j+1} s_{[1, j]} \i
s_{[j+2, a_{j+1}]} s_{j+1}=v_{j+1} s_{[j+2, a_{j+1}]} s_{[1, j+1]}
\i.}$$

Thus, we can prove by induction on $j$ that $x \i v_1 x=x_j \i v_j
s_{[j+1, a_j]} s_{[1, j]} \i x_j$ for $1 \le j \le i$. In
particular, $x \i v_1 x=s_{[i+1, n-1]} s_{[1, i]} \i$. The lemma
is proved. \qed

\subhead 3.7 \endsubhead By 3.4 and 3.5, there exists $t \in T$,
such that $(U \times U) (\dot v \dot x t, \dot x) \cdot h_J
\subset X$. Consider $K=\{a_j \mid b_{a_j}=0\}$. Then for any $j,
j' \in K$ with $j \neq j'$, we have that $|j-j'| \ge 2$ and
$<\a^\vee_j, \a_{j'}>=0$. Set $y=\prod_{j \in K} s_j$. Then $y$ is
well-defined. Note that $(1-y) y x \o^\vee_i, (1-y) v x \o^\vee_i
\in \sum_{j \in K} \bold R_{>0} \a^\vee_j$. By 3.3, $(\dot x \i
\dot y \dot v \dot x t, 1) \cdot h_J \in X$. Therefore, $(\dot x
\i \dot y \dot v \dot x t, 1) \cdot h_J \in X$. By 3.6, $x \i y v
x=x \i v_1 x=w^J$. Therefore, $Z_{J, w^J} \cap X \neq
\varnothing$. By 3.1, $G$ acts transitively on $Z_{J, w^J}$.
Therefore $Z_{J, w^J} \subset X$.

\subhead 3.8 \endsubhead In this subsection, we assume that $G$ is
of type $C_n$ and set $$\e=\cases 1, & \hbox{ if } 2 \mid i; \cr
0, & \hbox{ otherwise}. \cr \endcases$$

Set $v=s_{n-i+1} s_{n-i+3} \cdots s_{n-\e}$, $x_1=s_{[n-i, n-1]}
\i s_{[n-i-1, n-2]} \i \cdots s_{[1, i]} \i$ and $x_2=s_{[n+\e-1,
n]} \i s_{[n+\e-3, n]} \i \cdots s_{[n-i+2, n]} \i$. Set
$\l=\a^\vee_{n-i+1}+\a^\vee_{n-i+3}+\cdots+\a^\vee_{n-\e}$. Then
we have that $(v, \l)$ is admissible.

Now set $\l'=\sum_{j \in I} \min(i, j) \a^\vee_j \in \bold N
\o^\vee_i$. Set $x_{1, j}=s_{[j-i+1, j]} \i s_{[j-i, j-1]} \i
\cdots s_{[1, i]} \i$ for $i-1 \le j \le n-1$, s. Then we can show
by induction that $x_{1, j} \l'=\sum_{k=1}^i k \a^\vee_{j-i+1+k}+i
\sum_{l=j+2}^n \a^\vee_l$. In particular, $x_1
\o^\vee_i=\sum_{k=1}^i k \a^\vee_{n-i+k}$.

For $0 \le j \le (i+\e-1)/2$, set $x_{2, j}=s_{[n-i+2 j, n]} \i
s_{[n-i+2 j-2, n]} \i \cdots s_{[n-i+2, n]} \i$. Then we can show
by induction that $x_{2, j} x_1 \l'=\sum_{k=0}^{j-1}
\a^\vee_{n-i+1+2 k}+\sum_{l=1}^{i-2 j} l \a^\vee_{n-i+2 j+l}$. In
particular, we have that $x_2 x_1 \l'=\l$. Therefore, there exists
$t \in T$, such that $(U, U) (\dot v \dot x_2 \dot x_1 t, \dot x_2
\dot x_1) \cdot h_J \subset X$.

Now set $y_1=s_{n+\e-1} s_{n+\e-3} \cdots s_{n-i}$ and $y_2=s_{[1,
n-i-1]}$. For $1 \le j \le n-i-1$, set $\b_k=-(v x_2 x_1) \i
\a_k=-\a_{k+i}$. Thus by 3.3, there exists $u \in U_{\b_1}
U_{\b_2} \cdots U_{\b_{n-i}}$, such that $(\dot x_1 \i \dot x_2 \i
\dot y_1 \dot y_2 \dot v \dot x_2 \dot x_1 u t, 1) \cdot h_J \in
X$.

For $0 \le j \le (i+\e-1)/2$, set $$\eqalignno{v_{2, j}=s_{[1,
n-i]} (s_{n-i+2} s_{n-i+4} \cdots s_{n-i+2 j}) (s_{n-i+1}
s_{n-i+3} \cdots s_{n-i+2 j-1}) s_{[n-i+2 j+1, n]} \i.}$$

It is easy to see that $s_{[n-i+2 j, n]} v_{2, j} s_{[n-i+2 j, n]}
\i=v_{2, j-1}$. Therefore, we can show by induction that $x_2 \i
y_1 y_2 v x_2=x_{2, j} \i v_{2, j} x_{2, j}$ for $0 \le j \le
(i+\e-1)/2$. In particular, $x_2 \i y_1 y_2 v x_2=s_{[1, n-i]}
s_{[n-i+1, n]} \i$.

For $i-1 \le j \le n-1$, set $v_{1, j}=s_{[1, j-i+1]} s_{[j+2, n]}
s_{[j-i+2, j+1]} \i$. Then we have that $s_{[j-i+1, j]} v_{1, j}
s_{[j-i+1, j]} \i=v_{1, j-1}$. Therefore, we can show by induction
that $x_1 \i s_{[1, n-i]} s_{[n-i+1, n]} \i x_1=x_{1, j} \i v_{1,
j} x_{1, j}$ for $i-1 \le j \le n-1$. In particular, $x_2 \i y_1
y_2 v x_2=s_{[i+1, n]} s_{[1, i]} \i=w^J$.

Moreover, $w_{n-i-k+1} \i w^{-n+i+k+1} \b_k=w_{n-i-k+1} \i
(-\a_{n-1})=-\sum_{l=n-k}^n \a_l$. Since $n-k \in
J_{n-i-k-1}-J_{n-i-k}$, $U_{\b_k} \subset {}^{\dot w^{n-i-k-1}
\dot w_{n-i-k+1}} (U_{P^-_{J_{n-i-k}}} \cap U^-_{J_{n-i-k-1}})$.
By 3.1, $(\dot x_1 \i \dot x_2 \i \dot y_1 \dot y_2 \dot v \dot
x_2 \dot x_1 u t, 1) \cdot h_J \in Z_{J, w^J}$. Therefore, $Z_{J,
w^J} \subset X$.

${}$

For type $B_n$, we have the similar results.

\subhead 3.9 \endsubhead In subsection 3.9 and 3.10, we assume
that $G$ is of type $D_n$. In this subsection, assume that $i \le
n-2$.

If $2 \mid i$, set $v=s_{n-i} s_{n-i+2} \cdots s_{n-2}$,
$\l=\a^\vee_{n-i}+\a^\vee_{n-i+2}+\cdots+\a^\vee_{n-2}$ and
$x=(s_{[n-1, n]} \i s_{[n-3, n]} \i \cdots s_{[n-i+1, n]} \i)
(s_{[n-i-1, n-2]} \i s_{[n-i-2, n-3]} \i \cdots s_{[1, i]} \i)$.

If $2 \nmid i$, set $v=(s_{n-i} s_{n-i+2} \cdots s_{n-1}) s_n$,
$\l=\sum_{l=0}^{(i-3)/2} \a^\vee_{n-i+2
l}+1/2(\a^\vee_{n-1}+\a^\vee_n)$ and $x=(s_{[n-2, n]} \i s_{[n-4,
n]} \i \cdots s_{[n-i+1, n]} \i) (s_{[n-i-1, n-2]} \i s_{[n-i-2,
n-3]} \i \cdots s_{[1, i]} \i)$.

By the similar calculation to what we did for type $C_{n-1}$, we
have that in both cases $(v, \l)$ is admissible and $x \i
\l=\o^\vee_i$. Moreover, by the similar argument to what we did
for type $C_{n-1}$, we can show that $Z_{J, w^J} \subset X$.

\subhead 3.10 \endsubhead Assume that $i=n$. Set
$$\e=\cases 1, & \hbox{ if } 2 \mid [n/2]; \cr 0, & \hbox{
otherwise}. \cr
\endcases$$

If $2 \nmid n$, set $v=s_{n+\e-1} (s_1 s_3 \cdots s_{n-2})
s_{n-\e}$, $x=s_{n+\e-1} (s_{[n-3, n]} \i s_{[n-5, n]} \i \cdots
s_{[2, n]} \i) s_{n-1}$ and $\l=\frac {3}{2} \a^\vee_{n-\e}+\frac
{1}{2} \a^\vee_{a+\e-1}+\sum_{j=0}^{(n-3)/2} \a^\vee_{2 j+1}$.
Then $\l=2 x \o^\vee_n$ and $(v, \l)$ is admissible. Set $y=s_2
s_4 \cdots s_{n-3}$. Then $(\dot v \dot x t, \dot y \i \dot x)
\cdot h_J \in X$ for some $t \in T$. By 3.3, $(\dot x \i \dot y
\dot v \dot x t, 1) \cdot h_J \in X$. Since $x \i y v x=s_{n-1}
s_{[1, n-2]} \i s_n=w^J$, $Z_{J, w^J} \subset X$.

If $2 \mid n$, set $v=(s_1 s_3 \cdots s_{n-3}) s_{n-\e}$,
$\l=\a^\vee_{n-\e}+\sum_{j=0}^{n/2-2} \a^\vee_{1+2 j}$ and
$$x=\cases s_2 s_4, & \hbox{ if } n=4; \cr
s_{n-2} s_{n+\e-1} (s_{[n-4, n]} \i s_{[n-6, n]} \i \cdots s_{[2,
n]} \i) s_{n-1}, & \hbox{ otherwise}. \cr \endcases$$

Then $\l=2 x \o^\vee_n$ and $(v, \l)$ is admissible. Therefore,
there exists $t \in T$, such that $(U, U) (\dot v \dot x t, \dot
x) \cdot h_J \subset X$. Set $y_1=s_2 s_4 \cdots s_{n-2}$,
$y_2=s_{n+\e-1}$ and $\b=-(v x) \i \a_{n+\e-1}=-\a_{n/2}$. By 3.3,
there exists $u \in U_{\b}$ and $t \in T$, such that $(\dot x \i
\dot y_1 \dot y_2 \dot v \dot x u t, 1) \cdot h_J \in X$.

It is easy to see that $x \i y_1 y_2 v x=s_{n-1} s_{[1, n-2]} \i
s_n=w^J$ and $$w_2 \i \b=\cases -\sum_{l=1}^3 \a_l, & \hbox{ if }
n=4; \cr -\sum_{l=n/2-1}^{n-2} \a_l, & \hbox{ otherwise}. \cr
\endcases$$

Note that $J_0=I-\{n\}$ and $J_1=I-\{n-2, n\}$. Thus $U_{\b}
\subset {}^{w_2} (U_{P^-_{J_1}} \cap U^-_{J_0})$. By 3.1, $Z_{J,
w^J} \subset X$.

Similarly, $Z_{I-\{i-1\}, s_n s_{[1, n-2]} \i s_{n-1}} \subset X$.

\subhead 3.11. Type $G_2$ \endsubhead

Set $v=s_i$, $x=w^J$ and $\l=\a^\vee_i=x \o^\vee_i$. Then $(v,
\l)$ is admissible. Set $y=s_{3-i}$, then $(\dot x \i \dot y \dot
v \dot x t, 1) \cdot h_J \in X$ for some $t \in T$. Note that $x
\i y v x=w^J$. Therefore, $Z_{J, w^J} \subset X$.

\subhead 3.12. Type $F_4$ \endsubhead

If $i=1$, then set $v=s_2$, $x=s_1 s_4 w^2$ and $\l=\a^\vee_2=x
\o^\vee_1$. Thus $(v, \l)$ is admissible. Set $y_1=s_1 s_3$,
$y_2=s_4$ and $\b=-(v x) \i \a_4=-(\a_2+\a_3)$. Then there exists
$u \in U_{\b}$ and $t \in T$, such that $(\dot x \i \dot y_1 \dot
y_2 \dot v \dot x u t, 1) \cdot h_J \in X$. Note that $x \i y_1
y_2 v x=w^J$ and $w_2 \i \b=-(\a_2+2 \a_3+\a_4)$. By 3.1, $Z_{J,
w^J} \subset X$.

If $i=2$, then set $v=s_1 s_3$, $x=s_2 w^2$ and $\l=\a^\vee_1+
\a^\vee_3=x \o^\vee_2$. Thus $(v, \l)$ is admissible. Set $y=s_2
s_4$, then $(\dot x \i \dot y \dot v \dot x t, 1) \cdot h_J \in X$
for some $t \in T$. Note that $x \i y v x=w^J$. Thus $Z_{J, w^J}
\subset X$.

If $i=3$, then set $v=s_2 s_4$, $x=s_3 w^2$ and $\l=2
\a^\vee_2+\a^\vee_4=x \o^\vee_3$. Thus $(v, \l)$ is admissible.
Set $y=s_1 s_3$, then $(\dot x \i \dot y \dot v \dot x t, 1) \cdot
h_J \in X$ for some $t \in T$. Note that $x \i y v x=w^J$. Thus
$Z_{J, w^J} \subset X$.

If $i=4$, then set $v=s_3$, $x=s_1 s_4 w^2$ and $\l=\a^\vee_3=x
\o^\vee_1$. Thus $(v, \l)$ is admissible. Set $y_1=s_2 s_4$,
$y_2=s_1$ and $\b=-(v x) \i \a_1=-(\a_2+2 \a_3)$. Then there
exists $u \in U_{\b}$ and $t \in T$, such that $(\dot x \i \dot
y_1 \dot y_2 \dot v \dot x u t, 1) \cdot h_J \in X$. Note that $x
\i y_1 y_2 v x=w^J$ and $w_2 \i \b=-(\a_1+2 \a_2+2 \a_3)$. By 3.1,
$Z_{J, w^J} \subset X$.

\subhead 3.13. Type $E_6$ \endsubhead

If $i=1$, then set $v=s_1 s_5 s_3 s_6$, $x=s_1 s_4 s_3 s_1 s_6
w^J$ and $\l=\a^\vee_1+2 \a^\vee_3+\a^\vee_5+2 \a^\vee_6=3 x
\o^\vee_1$. Thus $(v, \l)$ is admissible. Set $y_1=s_4$, $y_2=s_2$
and $\b=-(v x) \i \a_2=-(\a_3+\a_4+\a_5)$. Then there exists $u
\in U_{\b}$ and $t \in T$, such that $(\dot x \i \dot y_1 \dot y_2
\dot v \dot x u t, 1) \cdot h_J \in X$. Note that $x \i y_1 y_2 v
x=w^J$ and $w_2 \i \b=-(\a_2+\a_3+2 \a_4+\a_5+\a_6)$. By 3.1,
$Z_{J, w^J} \subset X$.

Similarly, $Z_{I-\{6\}, s_2 s_1 s_3 s_4 s_5 s_6} \subset X$.

If $i=2$, then set $v=s_4$, $x=s_2 s_3 s_5 s_4 s_2 w^J$ and
$\l=\a^\vee_4=x \o^\vee_1$. Thus $(v, \l)$ is admissible. Set
$y_1=s_2 s_3 s_5$, $y_2=s_1 s_6$, $\b_1=-(v x) \i
\a_1=-(\a_4+\a_5)$ and $\b_2=-(v x) \i \a_6=-(\a_3+\a_4)$. Then
there exists $u \in U_{\b_1} U_{\b_2}$ and $t \in T$, such that
$(\dot x \i \dot y_1 \dot y_2 \dot v \dot x u t, 1) \cdot h_J \in
X$. Note that $x \i y_1 y_2 v x=w^J$, $w_2 \i \b_1=-\sum_{l=3}^6
\a_l$ and $w_2 \i \b_2=-(\a_1+\a_3+\a_4+\a_5)$. By 3.1, $Z_{J,
w^J} \subset X$.

If $i=3$, then set $v=s_3 s_6 s_1 s_4 s_5$, $x=s_2 s_3 s_4 s_1 s_3
w^J$ and $\l=2 \a^\vee_1+\a^\vee_3+3 \a^\vee_4+5
\a^\vee_5+\a^\vee_6=3 x \o^\vee_3$. Thus $(v, \l)$ is admissible.
Set $y=s_2$, then $(\dot x \i \dot y \dot v \dot x t, 1) \cdot h_J
\in X$ for some $t \in T$. Note that $x \i y v x=w^J$. Thus $Z_{J,
w^J} \subset X$.

Similarly, $Z_{I-\{5\}, s_2 s_1 s_3 s_4 s_6 s_5} \subset X$.

If $i=4$, then set $v=s_2 s_3 s_5$, $x=s_4 (w^J)^2$ and
$\l=\a^\vee_2+\a^\vee_3+5 \a^\vee_5=x \o^\vee_3$. Thus $(v, \l)$
is admissible. Set $y=s_1 s_4 s_6$, then $(\dot x \i \dot y \dot v
\dot x t, 1) \cdot h_J \in X$ for some $t \in T$. Note that $x \i
y v x=w^J$. Thus $Z_{J, w^J} \subset X$.

\subhead 3.14. Type $E_7$ \endsubhead

If $i=1$, then set $v=s_4$, $x=s_3 s_1 s_2 s_5 s_4 s_3 s_1 s_7
(w^J)^2$ and $\l=\a^\vee_4=x \o^\vee_1$. Thus $(v, \l)$ is
admissible. Set $y_1=s_3 s_2 s_5$ , $y_2=s_1 s_6 s_7$, $\b_1=-(v
x) \i \a_1=-\sum_{l=3}^6 \a_l$, $\b_2=-(v x) \i \a_6=-(\a_4+\a_5)$
and $\b_3=-(v x) \i \a_7=-(\a_2+\a_3+\a_4)$. Then there exists $u
\in U_{\b_3} U_{\b_2} U_{\b_1}$ and $t \in T$, such that $(\dot x
\i \dot y_1 \dot y_2 \dot v \dot x u t, 1) \cdot h_J \in X$. Note
that $x \i y_1 y_2 v x=w^J$, $w_2 \i \b_1=-\a_4-\sum_{l=2}^7
\a_l$, $w_2 \i \b_2=-\sum_{l=2}^6 \a_l$ and $w_3 \i (w^J) \i
\b_3=-(\a_2+\a_4+\a_5+\a_6)$. By 3.1, $Z_{J, w^J} \subset X$.

If $i=2$, then set $v=s_2 s_3 s_5 s_7$, $x=s_4 s_2 s_7 (w^J)^3$
and $\l=\a^\vee_2+2 \a^\vee_3+\a^\vee_5+\a^\vee_7=2 x \o^\vee_2$.
Thus $(v, \l)$ is admissible. Set $y=s_1 s_4 s_6$. Then $(\dot x
\i \dot y \dot v \dot x t, 1) \cdot h_J \in X$ for some $t \in T$.
Note that $x \i y v x=w^J$. Thus $Z_{J, w^J} \subset X$.

If $i=3$, then set $v=s_2 s_3 s_5$, $x=s_1 s_4 s_3 s_7 (w^J)^3$
and $\l=\a^\vee_2+\a^\vee_3+\a^\vee_5=x \o^\vee_3$. Thus $(v, \l)$
is admissible. Set $y_1=s_1 s_4 s_6$, $y_2=s_7$ and $\b=-(v x) \i
\a_7=-(\a_4+\a_5)$. Then there exists $u \in U_{\b_3} U_{\b_2}
U_{\b_1}$ and $t \in T$, such that $(\dot x \i \dot y_1 \dot y_2
\dot v \dot x u t, 1) \cdot h_J \in X$. Note that $x \i y_1 y_2 v
x=w^J$ and $w_2 \i \b=-(\a_2+\a_4+\a_5+\a_6)$. By 3.1, $Z_{J, w^J}
\subset X$.

If $i=4$, then set $v=s_1 s_4 s_6$, $x=s_2 s_3 s_5 s_4 (w^J)^3$
and $\l=\a^\vee_1+2 \a^\vee_4+\a^\vee_6=x \o^\vee_4$. Thus $(v,
\l)$ is admissible. Set $y=s_2 s_3 s_5 s_7$. Then $(\dot x \i \dot
y \dot v \dot x t, 1) \cdot h_J \in X$ for some $t \in T$. Note
that $x \i y v x=w^J$. Thus $Z_{J, w^J} \subset X$.

If $i=5$, then set $v=s_2 s_3 s_5 s_7$, $x=s_4 s_6 s_5 (w^J)^3$
and $\l=\a^\vee_2+2 \a^\vee_3+3 \a^\vee_5+\a^\vee_7=2 x
\o^\vee_5$. Thus $(v, \l)$ is admissible. Set $y=s_1 s_4 s_6$.
Then $(\dot x \i \dot y \dot v \dot x t, 1) \cdot h_J \in X$ for
some $t \in T$. Note that $x \i y v x=w^J$. Thus $Z_{J, w^J}
\subset X$.

If $i=6$, then set $v=s_4 s_6$, $x=s_1 s_5 s_7 s_6 (w^J)^3$ and
$\l=\a^\vee_4+\a^\vee_6=x \o^\vee_6$. Thus $(v, \l)$ is
admissible. Set $y_1=s_2 s_3 s_5 s_7$, $y_2=s_1$ and $\b=-(v x) \i
\a_1=-(\a_3+\a_4+\a_5)$. Then there exists $u \in U_{\b}$ and $t
\in T$, such that $(\dot x \i \dot y_1 \dot y_2 \dot v \dot x u t,
1) \cdot h_J \in X$. Note that $x \i y_1 y_2 v x=w^J$ and $w_2 \i
\b=-\a_4-\sum_{l=1}^5 \a_l$. By 3.1, $Z_{J, w^J} \subset X$.

If $i=7$, then set $v=s_2 s_5 s_7$, $x=s_6 s_7 s_4 s_5 s_6 s_7 s_1
(w^J)^2$ and $\l=\a^\vee_2+\a^\vee_5+\a^\vee_7=2 x \o^\vee_7$.
Thus $(v, \l)$ is admissible. Set $y_1=s_4 s_6$, $y_2=s_3 s_1$,
$\b_1=-(v x) \i \a_3=-(\a_3+\a_4+\a_5)$ and $\b_2=-(v x) \i
\a_1=-(\a_2+\a_4+\a_5+\a_6)$. Then there exists $u \in U_{\b_2}
U_{\b_1}$ and $t \in T$, such that $(\dot x \i \dot y_1 \dot y_2
\dot v \dot x u t, 1) \cdot h_J \in X$. Note that $x \i y_1 y_2 v
x=w^J$, $w_2 \i \b_1=-\a_4-\sum_{l=1}^6 \a_l$, $w_3 \i (w^J) \i
\b_2=-\a_4-\sum_{l=1}^5 \a_l$. By 3.1, $Z_{J, w^J} \subset X$.

\subhead 3.15. Type $E_8$ \endsubhead

If $i=1$, then set $v=s_4 s_6$, $x=s_3 s_1 s_2 s_5 s_4 s_3 s_1 s_8
(w^J)^5$ and $\l=\a^\vee_4+\a^\vee_6=x \o^\vee_1$. Thus $(v, \l)$
is admissible. Set $y_1=s_2 s_3 s_5 s_7$, $y_2=s_1 s_8$, $\b_1=-(v
x) \i \a_1=-\a_4-\sum_{l=2}^6 \a_l$ and $\b_2=-(v x) \i
\a_8=-\sum_{l=3}^7 \a_l$. Then there exists $u \in U_{\b_2}
U_{\b_1}$ and $t \in T$, such that $(\dot x \i \dot y_1 \dot y_2
\dot v \dot x u t, 1) \cdot h_J \in X$. Note that $x \i y_1 y_2 v
x=w^J$, $w_2 \i \b_1=-\a_4-\a_5-\sum_{l=2}^7 \a_l$ and $w_2 \i
\b_2=-\a_4-\sum_{l=2}^8 \a_l$. By 3.1, $Z_{J, w^J} \subset X$.

If $i=2$, then set $v=s_2 s_3 s_5 s_7$, $x=s_4 s_2 s_7 s_8
(w^J)^6$ and $\l=\a^\vee_2+\a^\vee_3+\a^\vee_5+\a^\vee_7=x
\o^\vee_2$. Thus $(v, \l)$ is admissible. Set $y=s_1 s_4 s_6 s_8$.
Then $(\dot x \i \dot y \dot v \dot x t, 1) \cdot h_J \in X$ for
some $t \in T$. Note that $x \i y v x=w^J$. Thus $Z_{J, w^J}
\subset X$.

If $i=3$, then set $v=s_2 s_3 s_5 s_7$, $x=s_1 s_4 s_3 s_7 s_8
(w^J)^6$ and $\l=\a^\vee_2+\a^\vee_3+2 \a^\vee_5+\a^\vee_7=x
\o^\vee_3$. Thus $(v, \l)$ is admissible. Set $y=s_1 s_4 s_6 s_8$.
Then $(\dot x \i \dot y \dot v \dot x t, 1) \cdot h_J \in X$ for
some $t \in T$. Note that $x \i y v x=w^J$. Thus $Z_{J, w^J}
\subset X$.

If $i=4$, then set $v=s_1 s_4 s_6 s_8$, $x=s_2 s_5 s_3 s_4 s_8
(w^J)^6$ and $\l=\a^\vee_1+3 \a^\vee_4+2 \a^\vee_6+\a^\vee_8=x
\o^\vee_4$. Thus $(v, \l)$ is admissible. Set $y=s_2 s_3 s_5 s_7$.
Then $(\dot x \i \dot y \dot v \dot x t, 1) \cdot h_J \in X$ for
some $t \in T$. Note that $x \i y v x=w^J$. Thus $Z_{J, w^J}
\subset X$.

If $i=5$, then set $v=s_2 s_3 s_5 s_7$, $x=s_4 s_6 s_5 (w^J)^6$
and $\l=\a^\vee_2+2 \a^\vee_3+2 \a^\vee_5+\a^\vee_7=x \o^\vee_5$.
Thus $(v, \l)$ is admissible. Set $y=s_1 s_4 s_6 s_8$. Then $(\dot
x \i \dot y \dot v \dot x t, 1) \cdot h_J \in X$ for some $t \in
T$. Note that $x \i y v x=w^J$. Thus $Z_{J, w^J} \subset X$.

If $i=6$, then set $v=s_1 s_4 s_6$, $x=s_1 s_5 s_7 s_6 (w^J)^6$
and $\l=\a^\vee_1+2 \a^\vee_4+\a^\vee_6=x \o^\vee_6$. Thus $(v,
\l)$ is admissible. Set $y_1=s_2 s_3 s_5 s_7$, $y_2=s_8$ and
$\b=-(v x) \i \a_8$. Then there exists $u \in U_{\b}$ and $t \in
T$, such that $(\dot x \i \dot y_1 \dot y_2 \dot v \dot x u t, 1)
\cdot h_J \in X$. Note that $x \i y_1 y_2 v x=w^J$ and $w_2 \i
\b=-\a_4-\sum_{l=1}^5 \a_l$. By 3.1, $Z_{J, w^J} \subset X$.

If $i=7$, then set $v=s_2 s_3 s_5$, $x=s_6 s_7 s_8 s_4 s_5 s_6 s_7
(w^J)^5$ and $\l=\a^\vee_2+\a^\vee_3+\a^\vee_5=x \o^\vee_7$. Thus
$(v, \l)$ is admissible. Set $y_1=s_1 s_4 s_6$, $y_2=s_7 s_8$,
$\b_1=-(v x) \i \a_7=-(\a_3+\a_4+\a_5)$ and $\b_2=-(v x) \i
\a_8=-(\a_2+\a_4+\a_5+\a_6)$. Then there exists $u \in U_{\b_2}
U_{\b_1}$ and $t \in T$, such that $(\dot x \i \dot y_1 \dot y_2
\dot v \dot x u t, 1) \cdot h_J \in X$. Note that $x \i y_1 y_2 v
x=w^J$, $w_2 \i \b_1=-\a_4-\sum_{l=1}^6 \a_l$ and $w_3 \i (w^J) \i
\b_2=-\a_4-\sum_{l=1}^5 \a_l$. By 3.1, $Z_{J, w^J} \subset X$.

If $i=8$, then set $v=s_4$, $x=s_1 s_5 s_6 s_7 s_8 (w^J)^5$ and
$\l=\a^\vee_4=x \o^\vee_8$. Thus $(v, \l)$ is admissible. Set
$y_1=s_5 s_2 s_3$, $y_2=s_1 s_6 s_7 s_8$, $\b_1=-(v x) \i
\a_1=-\a_4-\sum_{l=2}^7 \a_l$, $\b_2=-(v x) \i
\a_6=-(\a_3+\a_4+\a_5)$, $\b_3=-(v x) \i \a_7=w^J \b_2$ and
$\b_4=-(v x) \i \a_8=(w^J)^2 \b_2$. Then there exists $u \in
U_{\b_4} U_{\b_3} U_{\b_2} U_{\b_1}$ and $t \in T$, such that
$(\dot x \i \dot y_1 \dot y_2 \dot v \dot x u t, 1) \cdot h_J \in
X$. Note that $x \i y_1 y_2 v x=w^J$, $w_2 \i \b_1=-\sum_{l=3}^6
\a_l-\sum_{l=1}^7 \a_l$, $w_2 \i \b_2=-\a_4-\sum_{l=1}^7 \a_l$,
$w_3 \i (w^J) \i \b_3=-\a_4-\sum_{l=1}^6 \a_l$ and $w_4 \i
(w^J)^{-2} \b_4=-\a_4-\sum_{l=1}^5 \a_l$. By 3.1, $Z_{J, w^J}
\subset X$.

\head 4. The explicit description of $\bar{\cu}$ \endhead

\subhead 4.1 \endsubhead We assume that $G^1$ is a disconnected
algebraic group such that its identity component $G^0$ is
reductive. Following \cite{St, 9}, an element $g \in G^1$ is
called quasi-semisimple if $gBg\i=B,gTg\i=T$ for some Borel
subgroup $B$ of $G^0$ and some maximal torus $T$ of $B$. We have
the following properties.

(a) {\it if $g$ is semisimple, then it is quasi-semisimple. } See
\cite{St, 7.5, 7.6}.

(b) {\it Let $g \in G^1$ is a quasi-semisimple element and $T_1$
be a maximal torus of $Z_{G^0}(g)^0$, where $Z_{G^0}(g)^0$ is the
identity component of $\{x \in G^0 \mid x g=g x\}$. Then any
quasi-semisimple element in $g G^0$ is $G^0$-conjugate to some
element of $g T_1$.} See \cite{L4, 1.14}.

(c) {\it $g$ is quasi-semisimple if and only if the
$G^0$-conjugacy class of $g$ is closed in $G^1$.} See \cite{Spa,
1.15(f)} for the if-part, the only-if-part is due to Lusztig in an
unpublished note. His proof is as follows.

\proclaim{Proposition(Lusztig)} Let $g \in G^1$. Let $cl_{G^0} g$
be the $G^0$-conjugacy class of $g$. Assume that $cl_{G^0} g$ is
closed. Then $g$ is quasi-semisimple.
\endproclaim

Proof. The proof is due to Lusztig.

By \cite{St 7.2}, we can find a Borel subgroup $B$ such that $g B
g \i=B$. Let $cl_B g$ be the $B$-conjugacy class of $g$. Since
$cl_B g \subset cl_{G^0} g$ and $cl_{G^0} g$ is closed, we see
that the closure of $cl_B g$ is contained in $cl_{G^0} g$. By
\cite{Spa 1.15(e)}, the closure of $cl_B g$ contains a
quasi-semisimple element. Hence $cl_{G^0} g$ contains a
quasi-semisimple element. Hence $g$ is quasi-semisimple. \qed

\subhead 4.2 \endsubhead Let $\rho_i: G@>>>GL(V_i)$ be the
irreducible representation of $G$ with lowest weight $-\o_i$ and
$\bar{\rho}_i: \bar{G}@>>>P \bigl(\End(V_i) \bigr)$ be the
morphism induced from $\rho_i$ (see \cite{DS, 3.15}). Let $\cn$ be
the subvariety of $\bar{G}$ consisting of elements such that for
all $i \in I$, the images under $\bar{\rho}_i$ are represented by
nilpotent endomorphisms of $V_i$. We have the following result.

\proclaim{Theorem 4.3} We have that
$$\bar{\cu}-\cu=\cn=\bigsqcup_{J \subsetneqq I} \bigsqcup_{w \in
W^J, \supp(w)=I} Z_{J, w}.$$
\endproclaim

Proof. By 2.11 and the results in section 3, we have that
$$\bigsqcup_{J \subsetneqq I} \bigsqcup_{w \in W^J, \supp(w)=I}
Z_{J, w} \subset \bar{\cu}-\cu.$$

For $i \in I$, let $X_i$ be the subvariety of $P \bigl(\End(V_i)
\bigr)$ consisting of the elements that can be represented by
unipotent or nilpotent endomorphisms of $V_i$. Then $X_i$ is
closed in $P(\End(V_i))$. Thus, $\bar{\rho}_i(z) \in X_i$ for $z
\in \bar{\cu}$. Moreover, since $G$ is simple, for any $g \in
\bar{G}$, $\bar{\rho}_i(g)$ is represented by an automorphism of
$V_i$ if and only if $g \in G$. Thus if $z \in \bar{\cu}-\cu$,
then $\bar{\rho}_i(z)$ is represented by an nilpotent endomorphism
of $V_i$. Therefore $\bar{\cu}-\cu \subset \cn$.

Assume that $w \in W^J$ with $\supp(w) \neq I$ and $\cn \cap Z_{J,
w} \neq \varnothing$. Let $C$ be the closed $L_{J, w}$-stable
subvariety that corresponds to $\cn \cap Z_{J, w}$. We have seen
that $\dot w$ is a quasi-semisimple element of $N_G(L_{J, w})$.
Moreover, there exists a maximal torus $T_1$ in $Z_{L_{J,
w}}(w)^0$ such that $T_1 \subset T$. Since $C$ is an $L_{J,
w}$-stable nonempty closed subvariety of $C_{J, w}$, $\dot w t \in
C$ for some $t \in T_1$. Set $z=(\dot w t, 1) \cdot h_J$. Then $z
\in \cn$.

Since $\supp(w) \neq I$, there exists $i \in I$ with $i \notin
\supp(w)$. Then $-w \o_i=-\o_i$. Let $v$ be a lowest weight vector
in $V_i$. Assume that $\bar{\rho}_i(z)$ is represented by an
endomorphism $A$ of $V$. Then $A v \in k^* v$. Thus $z \notin
\cn$. That is a contradiction. Therefore $\cn \subset \bigsqcup_{J
\subsetneqq I} \bigsqcup_{w \in W^J, \supp(w)=I} Z_{J, w}$. The
theorem is proved. \qed

\subhead Remark \endsubhead Let $G=PGL_4(k)$ and $I=\{1, 2, 3\}$.
Then the theorem implies that $Z_{\{1, 3\}, s_2 s_1 s_3 s_2}
\subset \bar{\cu}$. By 2.5, we can see that $Z_{\{1, 3\}, s_2 s_1
s_3 s_2}$ contains infinitely many $G$-orbits. Therefore
$\bar{\cu}$ contains infinitely many $G$-orbits.

\proclaim{Corollary 4.4} Let $i \in I$ and $J=I-\{i\}$ and $w$ be
a Coxeter element of $W$ with $w \in W^J$. Then $\overline{Z_{J,
w}}=\bigsqcup_{K \subset J} \bigsqcup_{w' \in W^K, \supp(w')=I}
Z_{K, w'}$.
\endproclaim

Proof. Note that $Z_{J, w} \subset \bar{\cu} \cap (\bigsqcup_{K
\subset J} Z_K)$. Since $\bar{\cu}$ and $\bigsqcup_{K \subset J}
Z_K$ are closed, $\overline{Z_{J, w}} \subset \bar{\cu} \cap
(\bigsqcup_{K \subset J} Z_K)=\bigsqcup_{K \subset J}
\bigsqcup_{w' \in W^K, \supp(w')=I} Z_{K, w'}$. Therefore by 2.11,
$\overline{Z_{J, w}}=\bigsqcup_{K \subset J} \bigsqcup_{w' \in
W^K, \supp(w')=I} Z_{K, w'}$. \qed

\subhead 4.5 \endsubhead Let $\s: G@>>>T/W$ be the morphism which
sends $g \in G$ to the $W$-orbit in $T$ that contains an element
in the $G$-conjugacy class of the semisimple part $g_s$. The map
$\s$ is called Steinberg map. The fibers of $\s$ are called
Steinberg fibers. The unipotent variety is an example of Steinberg
fiber. Some other interesting examples are the regular semisimple
conjugacy classes of $G$.

Let $F$ be a fiber of $\s$. It is known that $F$ is a union of
finitely many $G$-conjugacy classes. Let $t$ be a representative
of $\s(F)$ in $T$, then $F=G_{diag} \cdot t U$ and
$\bar{F}=G_{diag} \cdot t \bar{U}$ (see \cite{Spr2, 1.4}). It is
easy to see that $t (\bar{U}-U) \subset \cn$. Thus
$\bar{F}-F=G_{diag} \cdot t (\bar{U}-U) \subset \cn$. Therefore,
if $(w, \l)$ is admissible and $x \i \cdot \l$ dominant, then
there exists some $t' \in T$ such that $(U \times U) (\dot w \dot
x t', \dot x) \cdot h_{I(x \i \l)} \subset t \bar{U}$. Thus by
2.11 and the results in section 3, $\bigsqcup_{J \subsetneqq I}
\bigsqcup_{w \in W^J, \supp(w)=I} Z_{J, w} \subset \bar{F}-F$.
Therefore, we have
$$\bar{F}-F=\cn=\bigsqcup_{J \subsetneqq I} \bigsqcup_{w \in W^J,
\supp(w)=I} Z_{J, w}.$$

Thus $\bar{F}-F$ is independent of the choice of the Steinberg
fiber $F$. As a consequence, in general, $\bar{F}$ contains
infinitely many $G$-orbits (answering a question that Springer
asked in \cite{Spr2}).

\subhead 4.6 \endsubhead For any variety $X$ that is defined over
the finite field $\bold F_q$, we write $|X|_q$ for the number of
$\bold F_q$-rational points in $X$.

If $G$ is defined and split over the finite field $\bold F_q$,
then for any $w \in W^J$, $|\tz_{J, w}|_q=|G|_q q^{-l(w)}$ (see
\cite{L4, 8.20}). Thus $$|Z_{J, w}|_q=|G|_q q^{-l(w)}
(q-1)^{-|I-J|}=(\sum_{u \in W} q^{l(u)}) (q-1)^{|J|} q^{l(w_0
w)}.$$

Set $L(w)=\{i \in I \mid w s_i<w\}$. Then $w \in W^J$ if and only
if $J \subset L(w_0 w)$. Moreover, if $w \neq 1$, then $L(w_0 w)
\neq I$. Therefore

$$\eqalignno{|\bar{\cu}-\cu|_q &=\sum_{J \neq I} \sum_{w \in W^J,
\supp(w)=I} |Z_{J, w}|_q=(\sum_{w \in W} q^{l(w)}) \sum_{J \neq I}
\sum_{w \in W^J, \supp(w)=I} (q-1)^{|J|} q^{l(w_0 w)} \cr
&=\sum_{w \in W} q^{l(w)} \sum_{\supp(w)=I} \sum_{J \subset L(w_0
w)} q^{l(w_0 w)} (q-1)^{|J|}\cr &=\sum_{w \in W} q^{l(w)}
\sum_{\supp(w)=I} q^{l(w_0 w)+|L(w_0 w)|}.}$$

\subhead Remark \endsubhead Note that $|\bar{G}|_q=\sum_{w \in W}
q^{l(w)} \sum_{w \in W} q^{l(w_0 w)+|L(w_0 w)|}$ (see \cite{DP,
7.7}). Our formula for $|\bar{\cu}-\cu|_q$ bears some resemblance
to the formula for $|\bar{G}|_q$.

\head Acknowledgements. \endhead I thank George Lusztig for
suggesting the problem and for many helpful discussions. I thank
T. A. Springer for his note \cite{Spr2}. I benefited a lot from
it. I also thank the referee for his/her help and suggestions.

\enddocument